\RequirePackage{fix-cm}
\documentclass[smallcondensed]{svjour3}     
\smartqed  
\usepackage{graphicx}
\usepackage{amsmath}
\usepackage{cite}
\usepackage{xcolor}
\usepackage[caption=false]{subfig}
\journalname{Acta Applicandae Math}
\begin{document}

\title{Novel Lagrangian Hierarchies, Generalized Variational ODE's and Families of Regular and Embedded Solitary Waves
}

\titlerunning{Novel Lagrangian Hierarchies\dots}        

\author{Ranses Alfonso-Rodriguez         \and
        S. Roy Choudhury 
}


\institute{Ranses Alfonso-Rodriguez \at
              Department of Mathematics, University of Central Florida, Orlando, FL, USA \\
              \email{ranses.alfonsorodriguez@ucf.edu}           
           \and
           S. Roy Choudhury \at
              Department of Mathematics, University of Central Florida, Orlando, FL, USA \\
              \email{sudipto.choudhury@ucf.edu}
}

\date{Received: date / Accepted: date}

\maketitle

\begin{abstract}
Hierarchies of Lagrangians of degree two, each only partly determined by the choice of  leading terms and with some coefficients remaining free, are considered. The free coefficients they contain satisfy the most general differential geometric criterion currently known for the existence of a Lagrangian and variational formulation, and derived by solution of the full inverse problem of the calculus of variations for scalar fourth-order ODEs respectively. However, our Lagrangians have significantly greater freedom since our existence conditions are for individual coefficients in the Lagrangian. In particular, the classes of Lagrangians derived here have four arbitrary or free functions, including allowing the leading coefficient in the resulting variational ODEs to be arbitrary, and with models based on the earlier general criteria for a variational representation being special cases. For different choices of leading coefficients, the resulting variational equations could also represent traveling waves of various nonlinear evolution equations, some of which recover known physical models. Families of regular and embedded solitary waves are derived for some of these generalized variational ODEs in appropriate parameter regimes, with the embedded solitons occurring only on isolated curves in the part of parameter space where they exist. Future work will involve higher order Lagrangians, the resulting equations of motion, and their solitary wave solutions.
\end{abstract}

\section{Introduction}
\label{sec:Intro}

Derivations and use of Lagrangians for higher-order differential equations have seen renewed interest in recent years. This is due to many recent applications including, but not limited to, higher-order field-theoretic models \cite{PUM, Man}, investigations of isochronous behaviors in a variety of nonlinear models \cite{TCG}, treatments of higher-order Painlev\'e equations \cite{CGK}, and variational treatments of embedded solitary waves of a variety of nonlinear wave equations \cite{VK}.

The problem of finding a Lagrangian for a given ordinary differential equation (ODE) is usually referred to as the inverse problem of the calculus of variations. The classical problem, restricting the Lagrangian to a function of the first derivatives of the dependent variable $u$ on $z$, has the form $u'' = F(z, u, u')$. The necessary and sufficient conditions for such an equation to be derivable from the Euler-Lagrange equation%
\begin{align} \label{EL1}
	\partial _{u}L& -  d(\partial _{u'}L)/dx = 0
\end{align}%
\noindent
are classical, and were first derived by Helmholtz \cite{Helm,Lop}, and also investigated by Darboux \cite{Dar}.

The inverse problem for a fourth-order ODE was treated by Fels \cite{Fels}. We shall consider his main result later, and compare it to those derived in this paper.

In this paper, we approach this inverse problem differently. By treating the coefficients of the highest derivative in the ODE as arbitrary functions of the dependent variable and its derivatives, we are able to derive Lagrangians for entire classes of differential equations. The resulting generalized families of Lagrangians may then be employed to investigate solutions of various members in the corresponding families of variational ODEs. In this paper, we employ them to derive regular and embedded solitary wave solutions of member equations of the ODE family.

The remainder of this paper is organized as follows. Section \ref{sec:2} considers a preliminary example where the Lagrangian follows easily by matching the terms in the ODE to those in the Euler-Lagrange equation. Less trivially, we then illustrate the procedure for constructing regular and embedded solitary wave solutions of the original family of ODEs using this Lagrangian.

Section \ref{sec:3} develops our new approach to the inverse problem for generalized classes of nonlinear fourth-order ODEs with leading-order coefficients that may be dependent on up to the second derivative of the dependent variable. By matching the ODE to the Euler-Lagrange variational equations, the general form of Lagrangian for ALL possible variational ODEs of this class is derived, together with the conditions on the ODE coefficients. Section \ref{sec:4} then employs the resulting Lagrangian to construct both regular and embedded solitary waves for a representative member of this most general variational ODE of this family. Section \ref{sec:5} shows that our conditions on the variational ODE coefficients satisfy the most general conditions currently known\cite{Fels} for variational fourth-order ODEs. However, our ODEs are  in fact significantly more general in that they leave several individual coefficients in the Lagrangian free, as well as allowing the resulting variational ODE to have arbitrary leading-order coefficients. Section \ref{sec:conclusions} briefly reviews the results and conclusions of the paper.


\section{An Example Lagrangian}
\label{sec:2}

To introduce some necessary techniques, let us first consider the traveling--wave equation given by
\begin{equation}\label{ode}
	u^{(4)}\left(d_2-d_3u'\right)-2d_3u''u'''+u''\left(d_1-2\mu u\right)-\mu {u'}^2=0.
\end{equation}

This equation is related to that obtained from equation (1) of Ilison et al\cite{Ili07}. Matching derivatives of various orders to the Euler-Lagrangian equation of form \eqref{EL1}, gives us the Lagrangian
\begin{equation}\label{lag}
	L=-\frac{d_1}{2}u_z^2-\frac{\mu}{2}u^2u_{zz}+\frac{1}{2}\left(d_2-d_3u_z\right)u_{zz}^2.
\end{equation}

In the next two sub-sections we introduce the techniques for finding both regular and embedded solitons for the equation \eqref{ode} using its Lagrangian \eqref{lag} \cite{Let09,Smi09,VPh}.

\subsection{Regular Solitons}
\label{subsec:2.1}

In this section, we shall use the well--known Rayleigh--Ritz method for constructing regular solitary waves with exponentially decaying tails of Equation \eqref{ode}, as widely employed in many areas of applications. The extension to embedded solitary waves in the next sub-section is much more recent, and also less widely known.

As usual, the localized regular solitary wave solutions will be found using a Gaussian trial function \eqref{gtf}. This is standard even for simpler nonlinear PDEs where exact solutions may be known and have the usual $sech$ or $sech^2$ functional forms. We use the form:

\begin{equation}\label{gtf}
	u=A\exp\left(-\frac{z^2}{\rho^2}\right).
\end{equation}

Then, after substituting the trial function into the Lagrangian and integrating over the entire z axis one gets the `averaged action'  \eqref{action1}:

\begin{equation}\label{action1}
	S=\dfrac{A^2\left(27\sqrt{2}d_2+\left(-9\sqrt{2}d_1+8\sqrt{3}A\mu\right)\right)\sqrt{\pi}}{36\rho^3}.
\end{equation}

At this point the goal is to optimize the trial function by varying the action with respect to the amplitude, $A$, and the width, $\rho$. This determines the optimal parameters for the trial function or solitary wave solution, but within the particular functional form chosen for the ansatz. The resulting variational Euler--Lagrange equations, by varying $A$ and $\rho$ respectively, are the system of algebraic equations:

\begin{subequations}
	\begin{align}
		\frac{2 \sqrt{\frac{\pi }{3}} A^2 \mu}{3 \rho}+\frac{\sqrt{\pi } A \left(\rho^2 \left(8 \sqrt{3} A \mu-9 \sqrt{2} d_1\right)+27 \sqrt{2} d_2\right)}{18 \rho^3} & =0, \\
		\frac{\sqrt{\pi } A^2 \left(8 \sqrt{3} A \mu-9 \sqrt{2} d_1\right)}{18 \rho^2}-\frac{\sqrt{\pi } A^2 \left(\rho^2 \left(8 \sqrt{3} A \mu-9 \sqrt{2} d_1\right)+27 \sqrt{2} d_2\right)}{12 \rho^4} & =0.
	\end{align}
\end{subequations}

\noindent
The nontrivial solution to these equations is
\begin{subequations}
	\begin{align}
		\label{solA}A= &\dfrac{3\sqrt{\frac{3}{2}}d_1}{7\mu},\\
		\label{solrho}\rho^2= & \frac{21d_2}{d_1}.
	\end{align}
\end{subequations}

\begin{figure*}[ht]
   \subfloat[The regular soliton plotted for $d_1=1$.]{\label{Fig1a}
      \includegraphics[width=.47\textwidth]{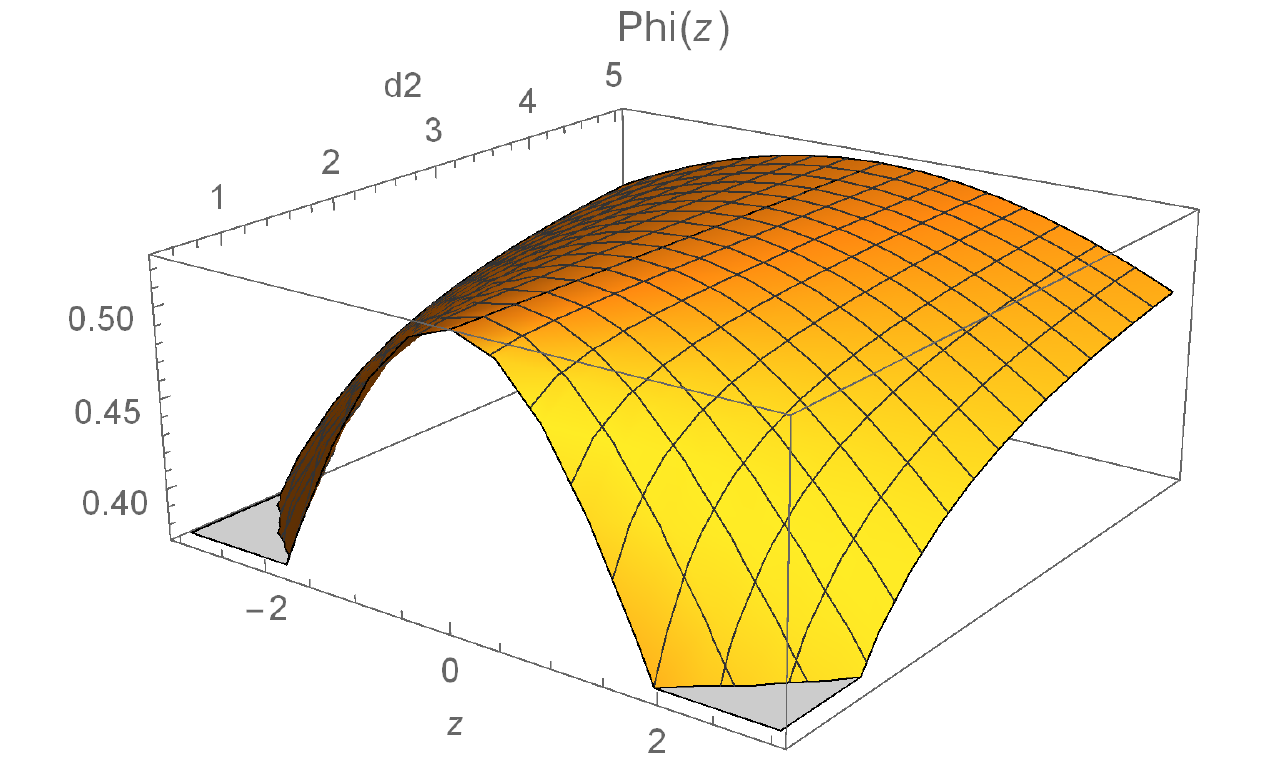}}
~
   \subfloat[The regular soliton plotted for $d_2=1$.]{\label{Fig1b}
     \includegraphics[width=.47\textwidth]{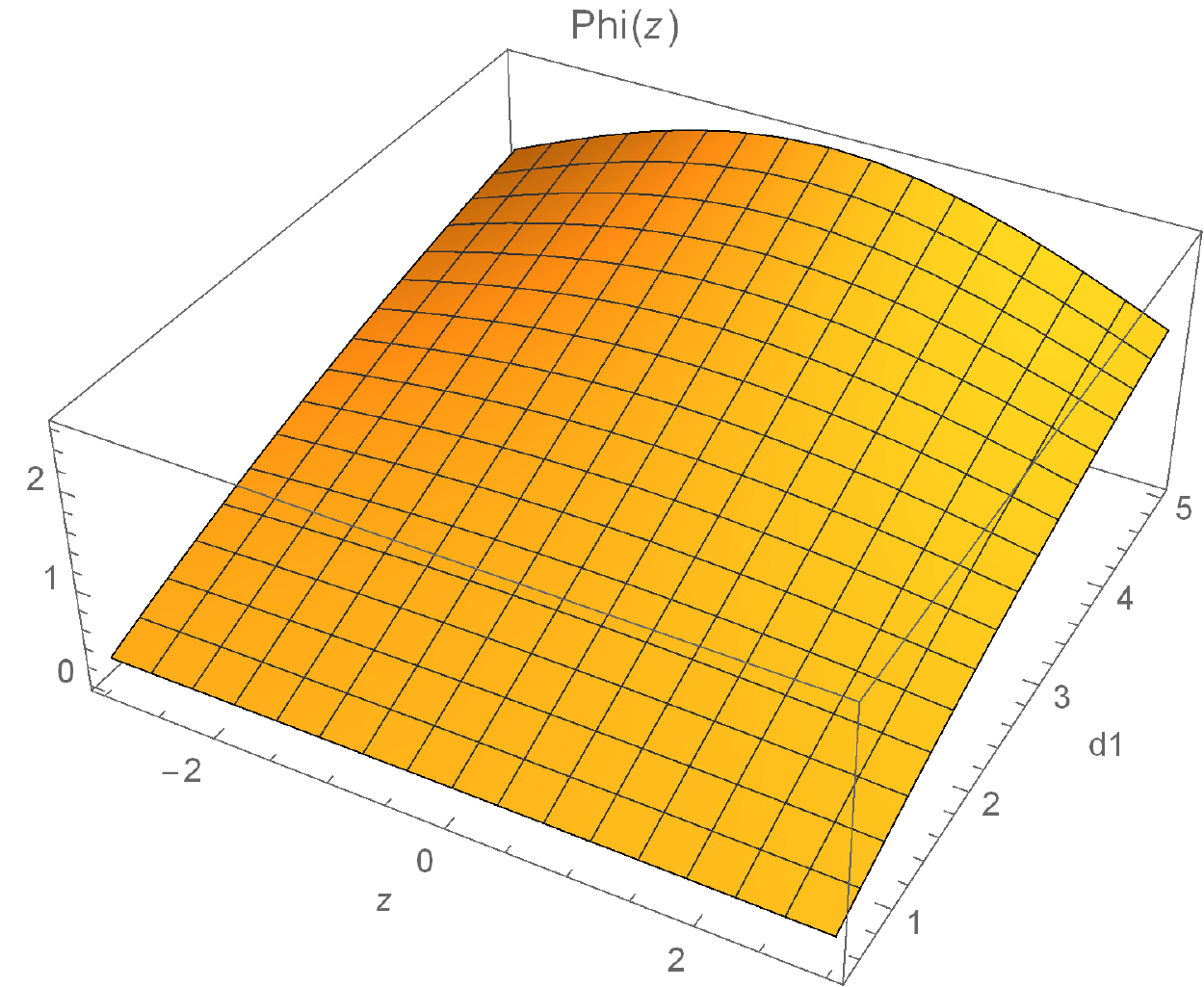}}

   \caption{Regular soliton for $z$ as one of the parameters varies and the other has value 1.}\label{Fig1}
\end{figure*}

The optimized variational soliton for the regular solitary waves of the traveling--wave equation \eqref{ode} is given by the trial function \eqref{gtf} with the above $A$ and $\rho$. Figures \ref{Fig1a} and \ref{Fig1b} respectively show the resulting regular solitary wave solution for $\mu=d_3=1$, and various positive values of either $d_2$ or $d_1$ (with the other one set to the value one).  



\begin{figure*}[ht]
   \subfloat[Direct accuracy analysis for $d_2=5$.]{\label{Fig2a}
      \includegraphics[width=.47\textwidth]{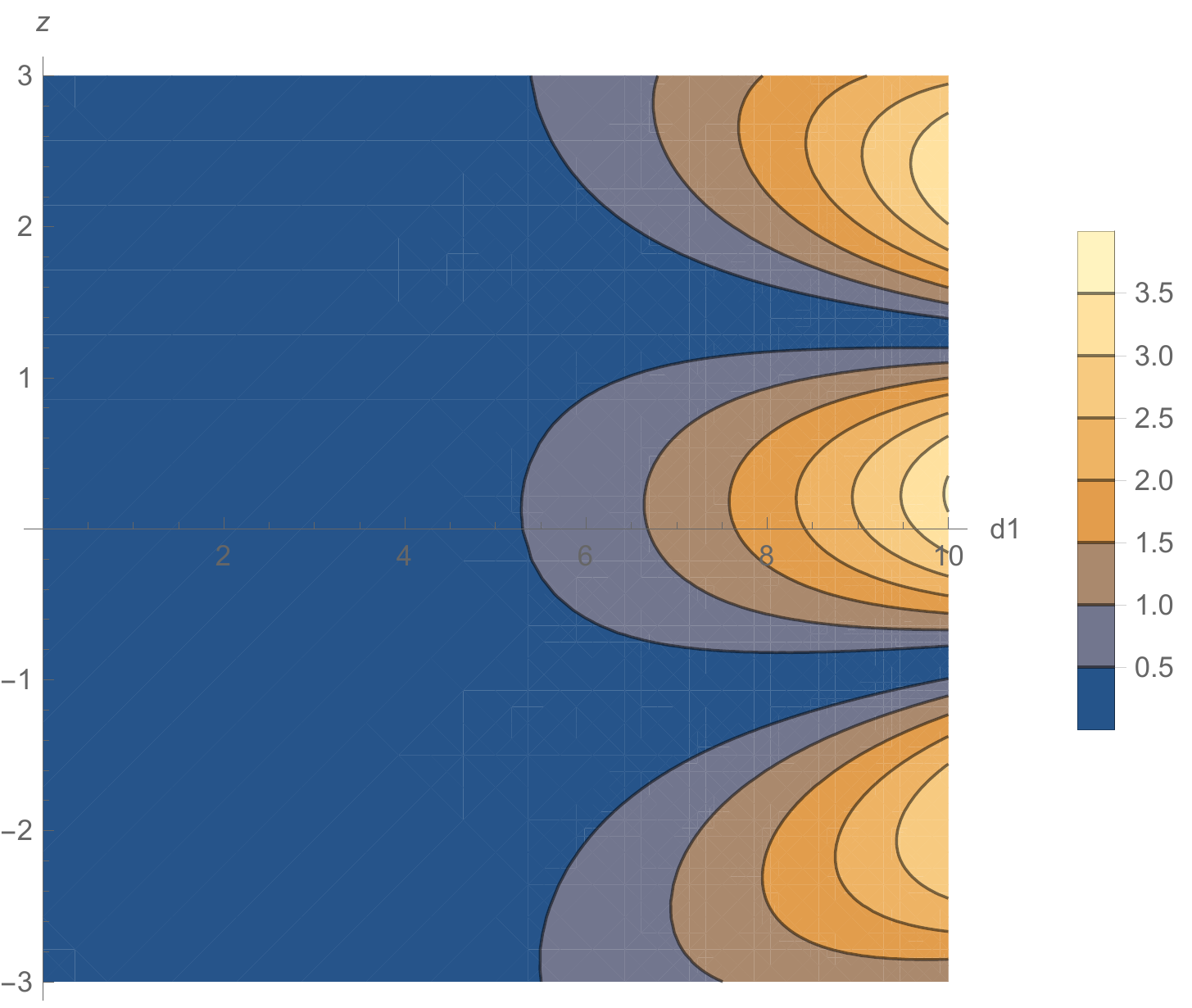}}
~
   \subfloat[Direct accuracy analysis for $d_2=1$.]{\label{Fig2b}
     \includegraphics[width=.47\textwidth]{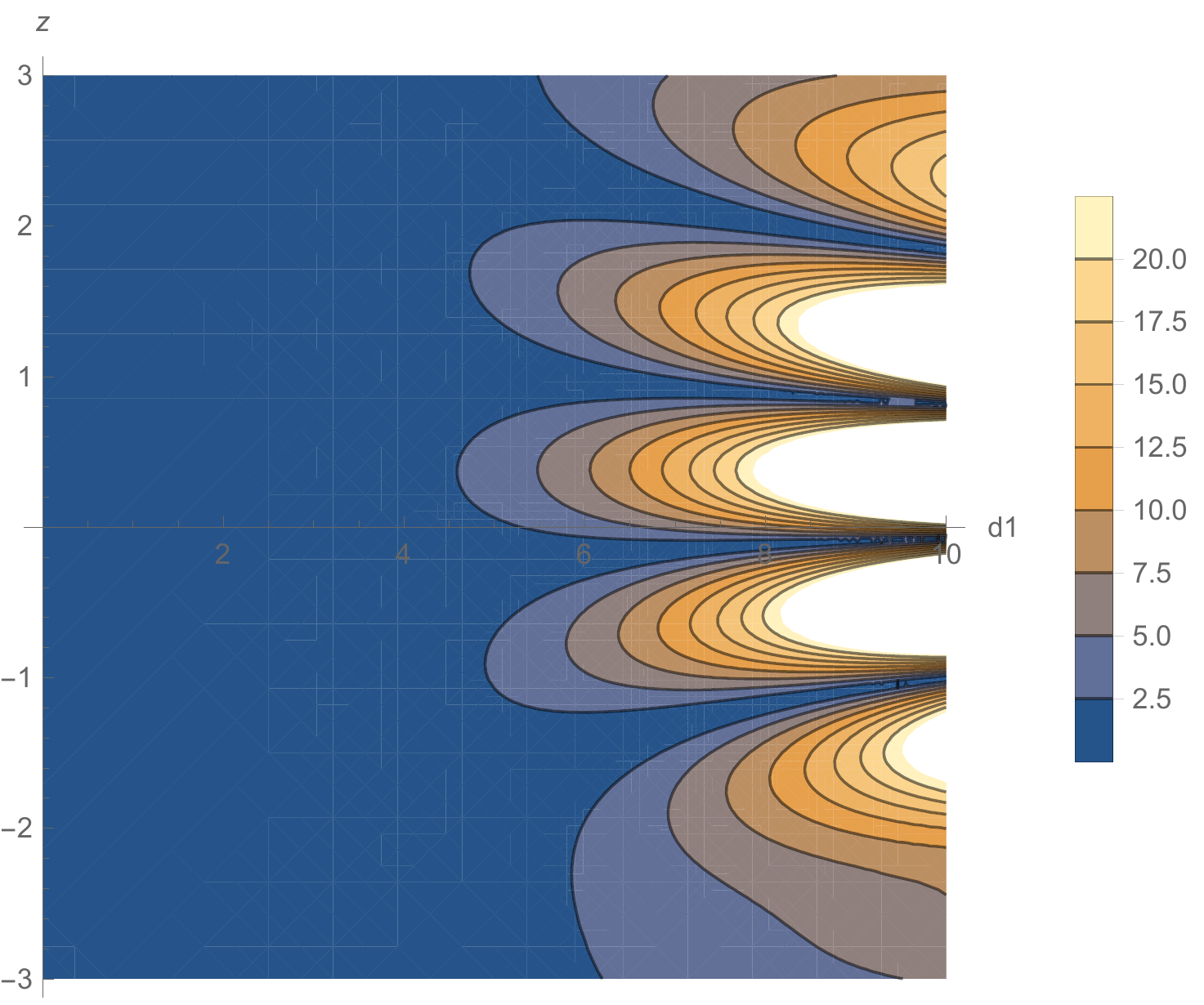}}

   \caption{Direct accuracy analysis of the regular soliton.}\label{Fig2}
\end{figure*}
\begin{figure}[htp]
	\centering
	\includegraphics[width=.5\textwidth]{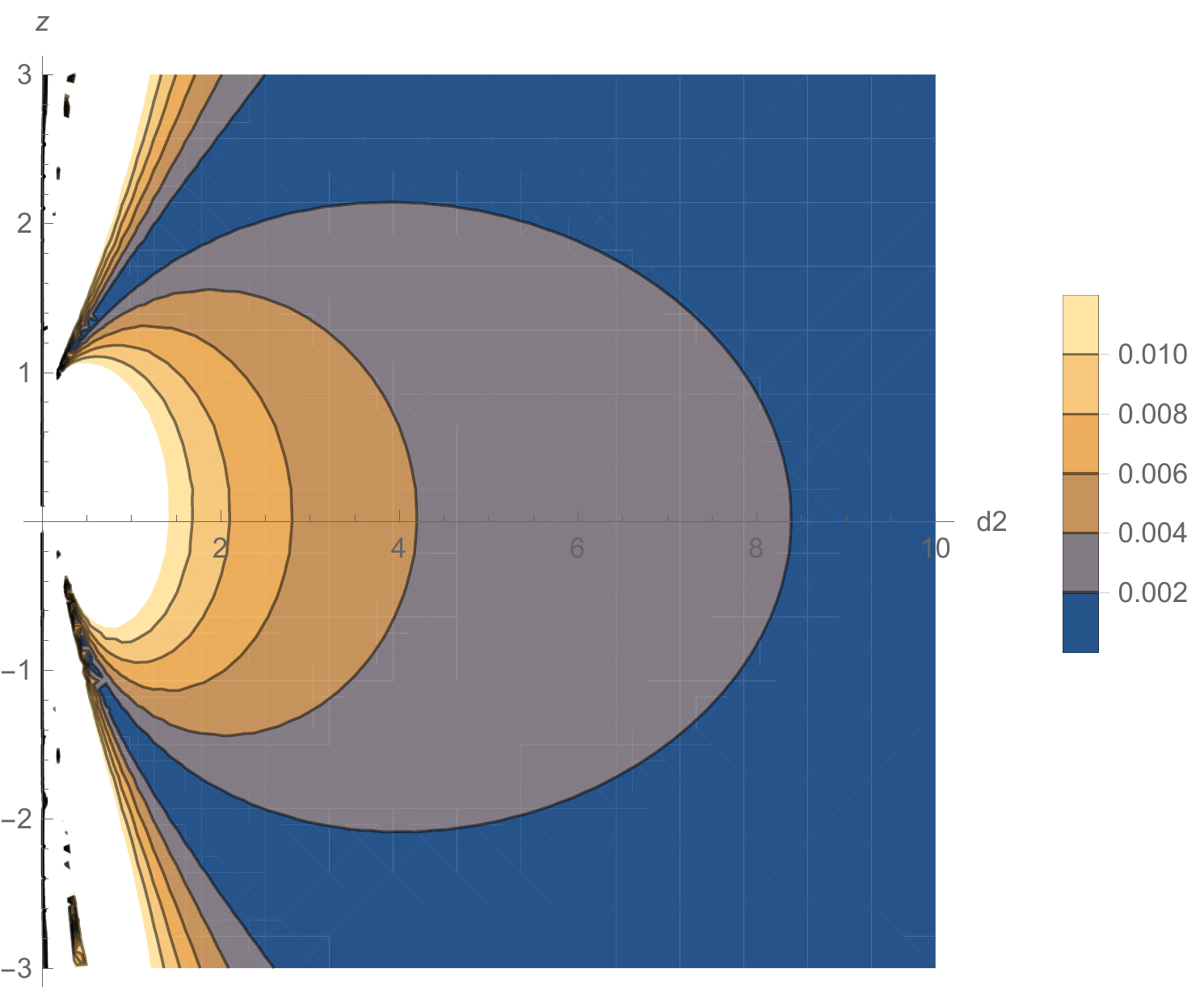}
	\caption{Direct accuracy analysis for $d_1=1$.}\label{Fig3}
\end{figure}

Figures \ref{Fig2} and \ref{Fig3} show a direct analysis of the accuracy of the variational regular solitary waves obtained above. It is possible to do this direct accuracy analysis since our variational solution for the regular solitary waves given by \eqref{gtf}, \eqref{solA}, and \eqref{solrho} is analytical. Substituting this variational solution \eqref{gtf} (with \eqref{solA} and \eqref{solrho}) into the traveling--wave \eqref{ode}, the deviation of the left-hand side of \eqref{ode} from zero gives a direct measure of the goodness of the variational solution. 

Figures \ref{Fig2a} and \ref{Fig2b} show the absolute value of this left-hand side for $\mu=1$, $d_3=1$ and $d_2=5$ and $d_2=1$ respectively for $z$ from -3 to 3 and $d_1$ from 0 to 10. It is clear that, for these fixed values, our variational solution is increasingly accurate when $d_1$ values are small.

Figure \ref{Fig3} displays the absolute value of this left-hand side of \eqref{ode} in a somewhat different fashion, in this case for $d_1 = 1$. Over the entire ranges of $z$ from -3 to 3 and $d_2$ from 0 to 10, it remains small, thus showing that our variational solution is accurate for all $z$ with small $d_2$ values. At larger $d_2$ values, the deviation from zero becomes even smaller, and the variational solution even more accurate.

\subsection{Embedded Solitons}
\label{subsec:2.2}

In the variational approach to embedded solitary waves, the tail of a delocalized soliton is modeled by 
\begin{equation}\label{destail}u_{\rm{tail}}=\alpha \cos(\kappa(c)z),\end{equation}
where the cosine ensures an even solution, and the arbitrary function $\kappa(c)$ will be sought so it ensures the integrability of the action \cite{VPh}.

Our ansatz for the embedded soliton uses a second order exponential core model plus the above tail model, i.e.,
\begin{equation}\label{aes}u=A \exp\left(-\frac{z^2}{\rho^2}\right)+u_{\rm{tail}}.\end{equation}

\noindent
Substituting this trial function into the Lagrangian \eqref{lag} and reducing the higher powers of the trigonometric functions yields an equation with trigonometric functions of the double and triple angles, as well as terms linear in $z$. The former would make averaging of the Lagrangian divergent. However, it is possibly to rigorously establish that such terms may be averaged out, so we shall set them to zero {\it a priori} \cite{Smi09, VPh}. 

The terms linear in $z$ would also cause the Lagrangian to be non-integrable. Hence, we set $\kappa(c)=\pm\sqrt{\frac{d_1}{d_2}}$, to force these linear terms to equal zero \cite{VPh}. This wavenumber of the tail oscillation \eqref{destail} may also be directly derived by substituting  \eqref{destail} into the linearized version of \eqref{ode} (for the small amplitude tail oscillations).

These two steps are not part of the traditional Rayleigh-Ritz method used for the construction of regular solitary waves and are used only for the variational approximation of embedded solitary waves.

Next, we can integrate the rest of the equation to obtain the action
\begin{multline}\label{action2}
	\frac{\left(\sqrt{\pi } A e^{-\frac{d_1 \rho^2}{8 d_2}}\right)}{72 d_2 \rho^3}\left(9 \sqrt{2} A \alpha \mu \rho^2 \left(3 d_1 \rho^2+4 d_2\right)\right.\\ \left. +2 e^{\frac{d_1 \rho^2}{8 d_2}} \left(8 \sqrt{3} A^2 d_2 \mu \rho^2+9 \sqrt{2} A d_2 \left(3 d_2-d_1 \rho^2\right)+18 \alpha^2 d_1 \mu \rho^4\right) \right)=0
\end{multline}

\noindent
As for the regular solitary waves, the action is now varied with respect to the amplitude ($A$), the width ($\rho$), and this time also the small amplitude ($\alpha$) of the oscillating tail to give the following system of equations:
\begin{subequations}\label{el}
    \begin{align} 
        0 & =8 \sqrt{3} A d_2 \mu \rho^2+6 \alpha^2 d_1 \mu \rho^4 \notag \\ \label{el1}
        & \qquad \qquad +3 \sqrt{2} \left(2 d_2 \rho^2 \left(2 \alpha \mu e^{-\frac{d_1 \rho^2}{8 d_2}}-d_1\right)+3 \alpha d_1 \mu \rho^4 e^{-\frac{d_1 \rho^2}{8 d_2}}+6 d_2^2\right) \\ 
        0 & =64 \sqrt{3} A^2 d_2^2 \mu \rho^2 e^{\frac{d_1 \rho^2}{8 d_2}}-144 \alpha^2 d_1 d_2 \mu \rho^4 e^{\frac{d_1 \rho^2}{8 d_1}} \notag \\ \label{el2}
        & \qquad \qquad +9 \sqrt{2} A \left(3 \alpha d_1^2 \mu \rho^6-8 d_2^2 \rho^2 \left(d_1 e^{\frac{d_1 \rho^2}{8 d_2}}-2 \alpha \mu\right)-8 \alpha d_1 d_2 \mu \rho^4+72 d_2^3 e^{\frac{d_1 \rho^2}{8 d_2}}\right) \\ \label{el3}
        0 & =\sqrt{2} A \left(3 d_1 \rho^2+4 d_1\right)+8 \alpha d_1 \rho^2 e^{\frac{d_1 \rho^2}{8 d_1}}
\end{align}
\end{subequations}

For embedded solitary waves, which live in a sea of delocalized solitary waves with small tails, the amplitude of the tail is strictly zero \cite{Let09, VPh}. Hence, we set $\alpha=0$ in the above equations in order to recover such embedded solitary waves, yielding
\begin{subequations}
	\begin{align}
		8 \sqrt{3} A d_2 \mu \rho^2+3 \sqrt{2} \left(6 d_2^2-2 d_1 d_2 \rho^2\right) & =0 \\
		64 \sqrt{3} A d_2^2 \mu \rho^2+9 \sqrt{2} \left(72 d_2^3-8 d_1 d_2^2 \rho^2\right) & =0\\
		\sqrt{2} A \left(3 d_1 \rho^2+4 d_2\right) & =0
	\end{align}
\end{subequations}

A nontrivial analytical solution to these equations is found to be
\begin{subequations}
	\begin{align}
		\label{sol2A}A= & \frac{13 \sqrt{\frac{3}{2}} d_1}{8 \mu} \\
		\label{sol2rho}\rho^2= & - \frac{4d_2}{3d_1}
	\end{align}
\end{subequations}

Note that using the model $exp(-z^2/\rho^2)$ for the core of the embedded solitary wave, with the tail \eqref{destail} added on results in an imaginary $\rho$ unless $d_1$ and $d_2$ have opposite signs, which would have made the argument to the exponential positive, and there would have been no decaying wave corresponding to a genuine embedded solitary wave.

\begin{figure*}[ht]
   \subfloat[The embedded soliton plotted for $\alpha=0$.]{\label{Fig4a}
      \includegraphics[width=.47\textwidth]{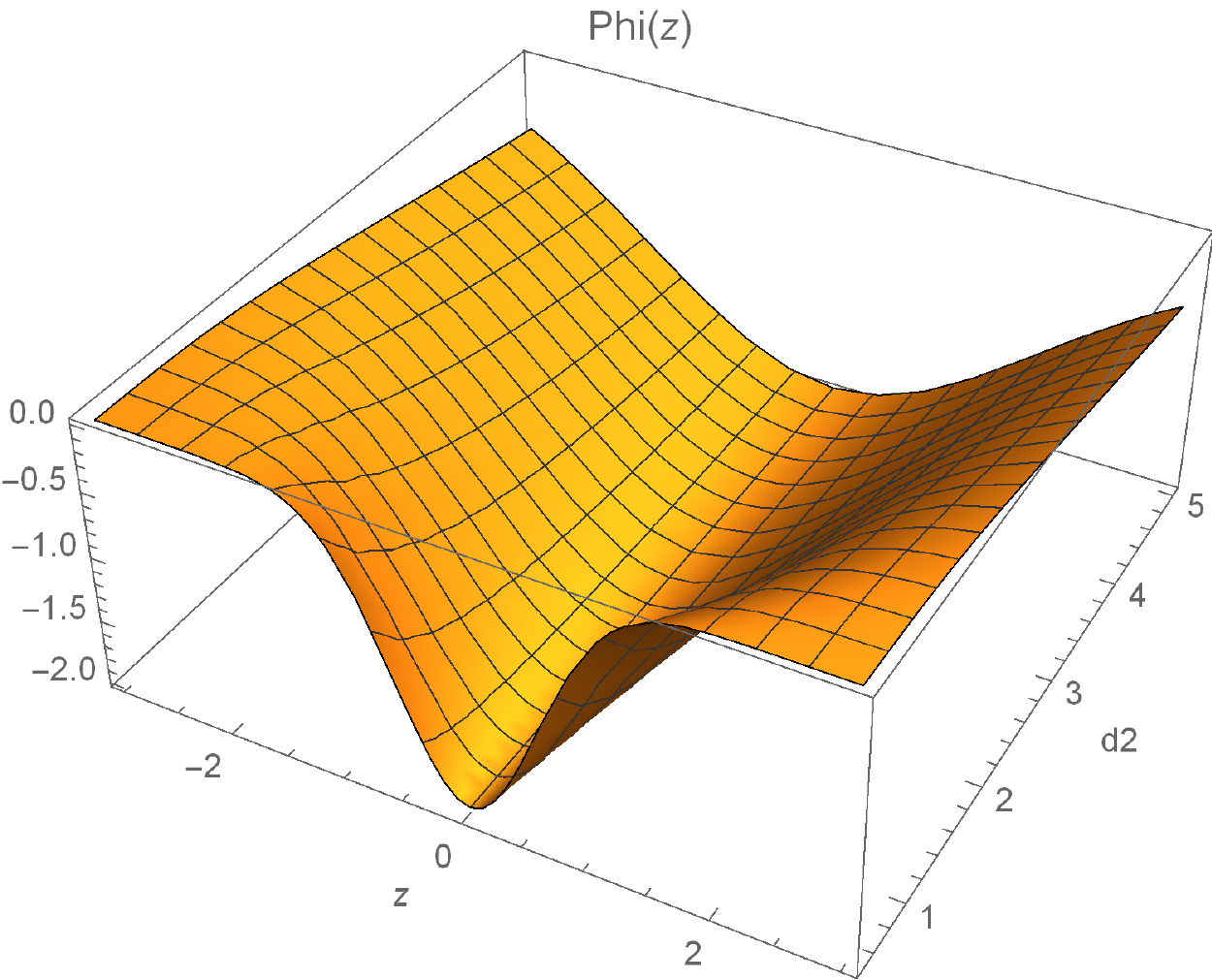}}
~
   \subfloat[The embedded soliton plotted for $\alpha=1$.]{\label{Fig4b}
     \includegraphics[width=.47\textwidth]{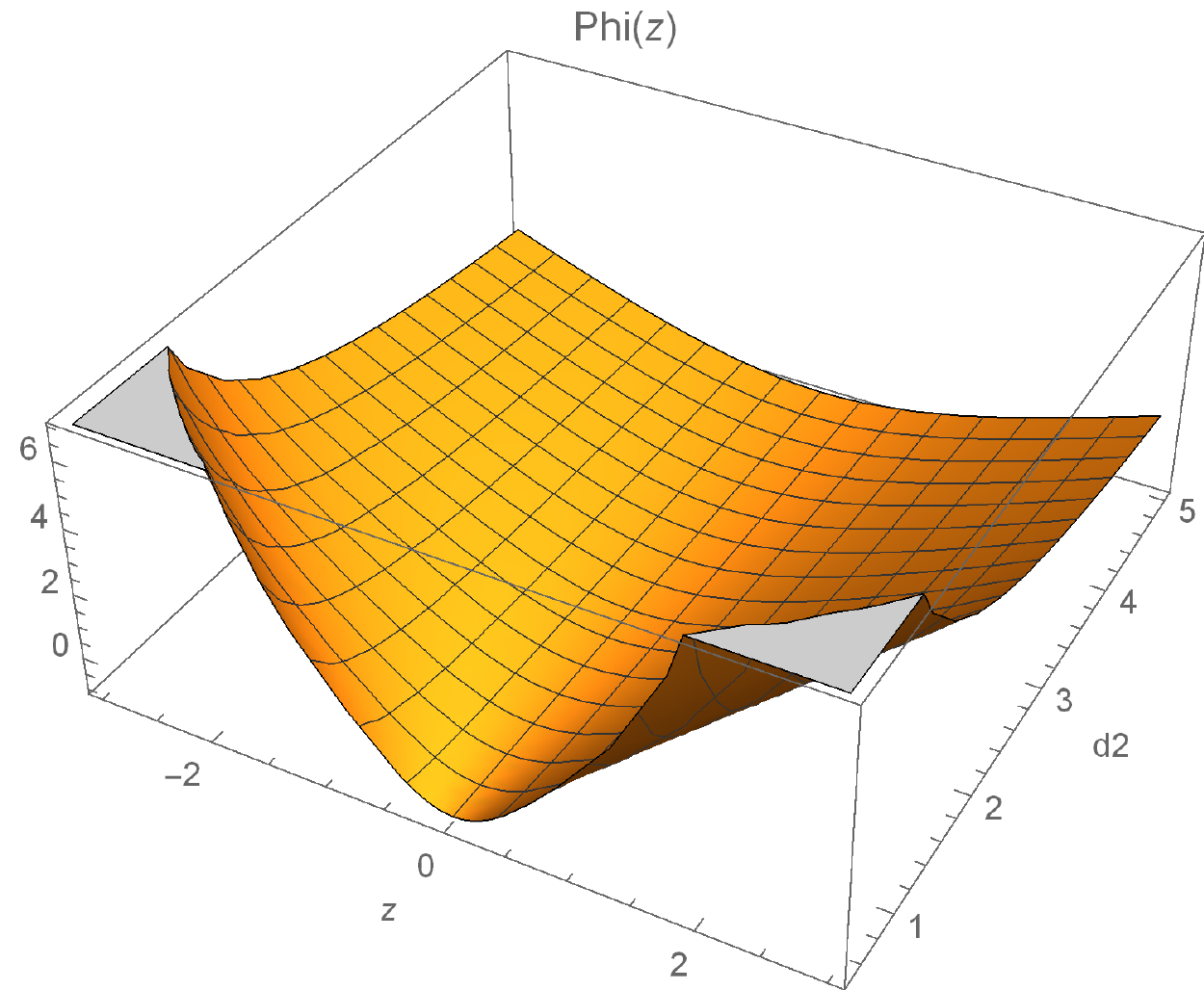}}

   \caption{Embedded soliton for $z$ and different values of the parameters.}\label{Fig4}
\end{figure*}

In Figures \ref{Fig4a} and \ref{Fig4b} the embedded soliton \eqref{aes}, with \eqref{sol2A}-\eqref{sol2rho}, is presented for $d_1=-1$, $\mu=1$ and $\alpha=0$ and $\alpha=1$ respectively, with $z$ varying from -3 to 3 and $d_2$ from 0.5 to 5.



\begin{figure*}[ht]
   \subfloat[Direct accuracy analysis of the embedded soliton for $d_2=-1$.]{\label{Fig5a}
      \includegraphics[width=.47\textwidth]{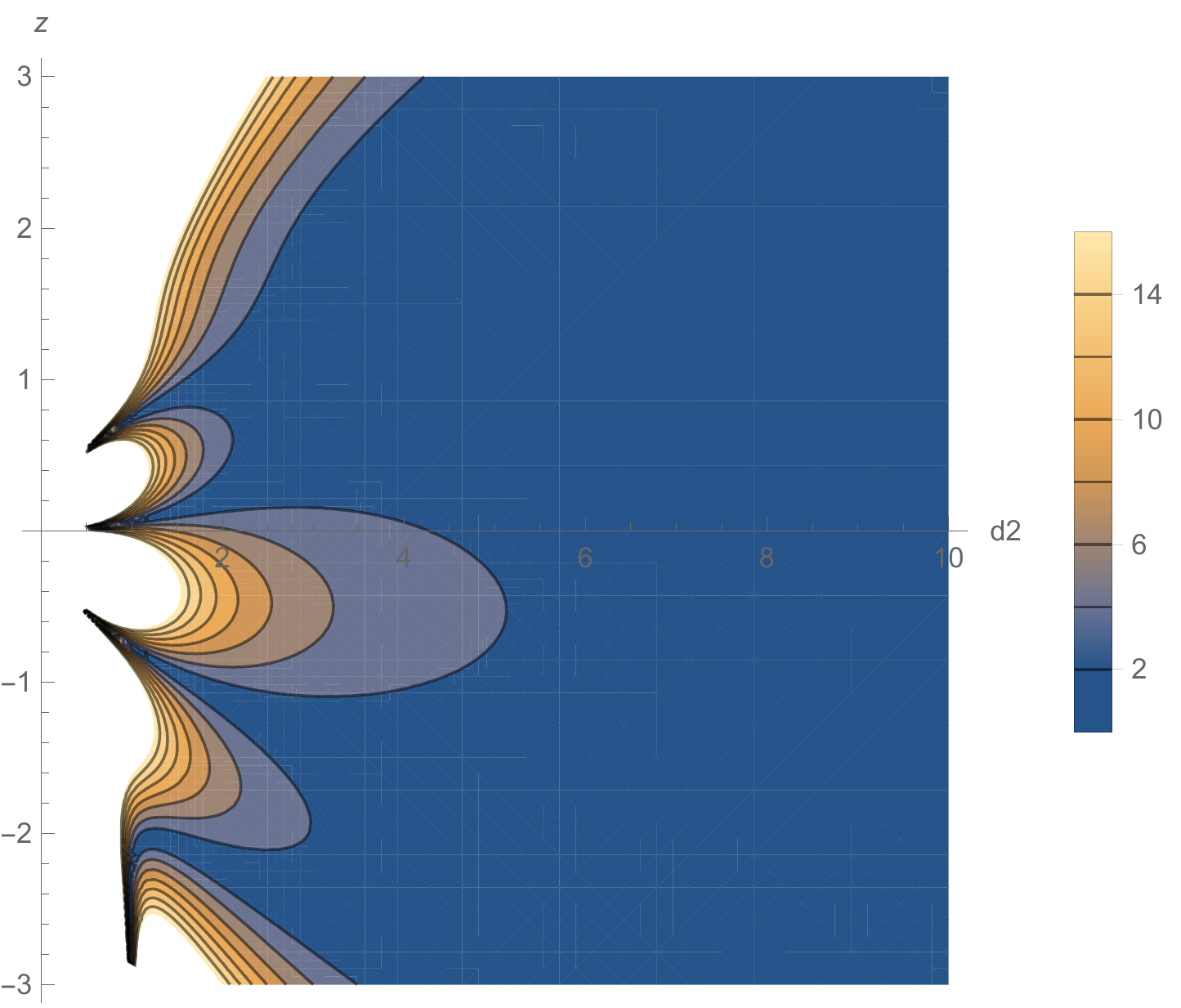}}
~
   \subfloat[Direct accuracy analysis of the embedded soliton for $d_2=-3$.]{\label{Fig5b}
     \includegraphics[width=.47\textwidth]{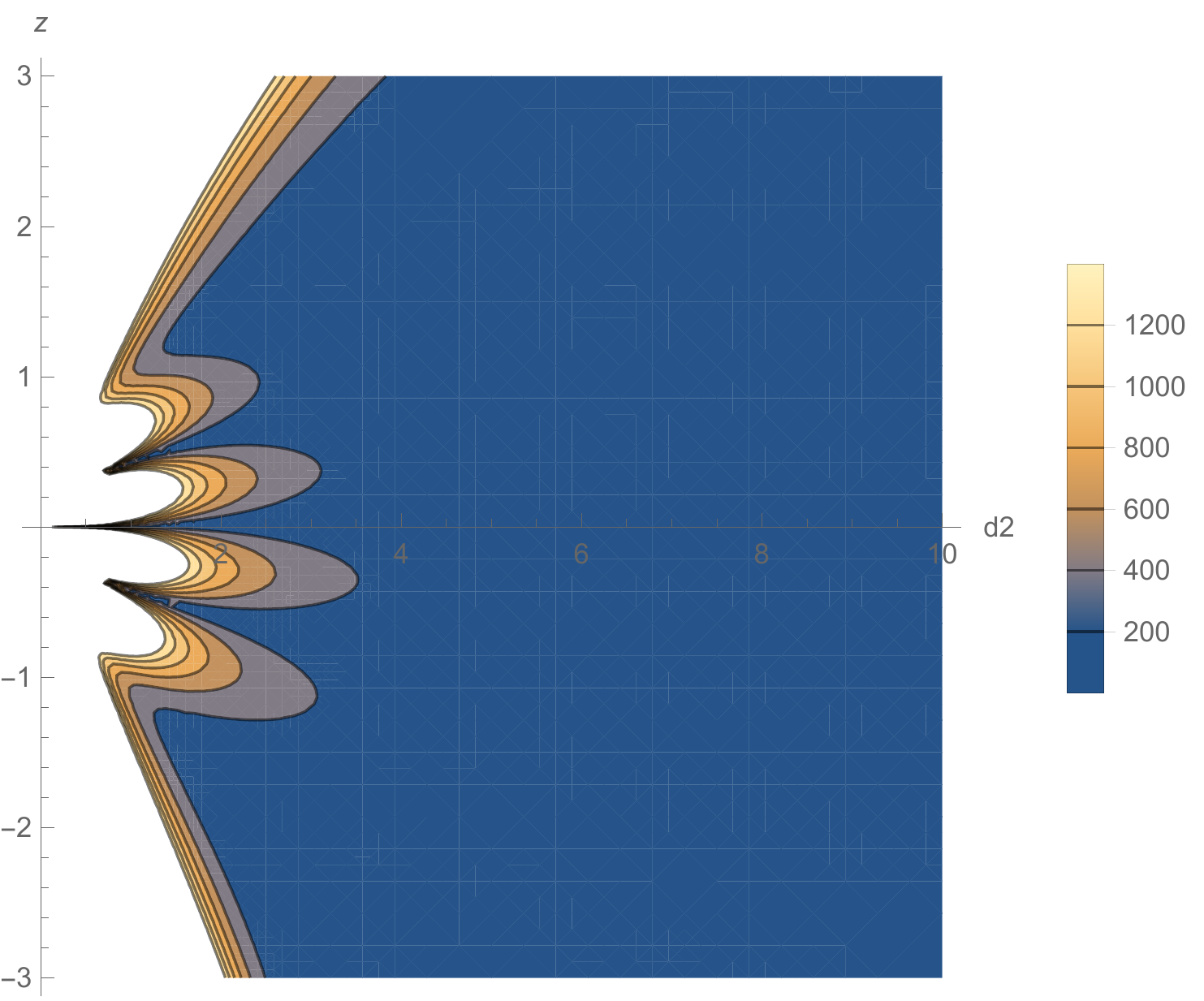}}
\\
   \centering
   \subfloat[Direct accuracy analysis of the embedded soliton for $d_2=-1/3$.]{\label{Fig5c}
     \includegraphics[width=.47\textwidth]{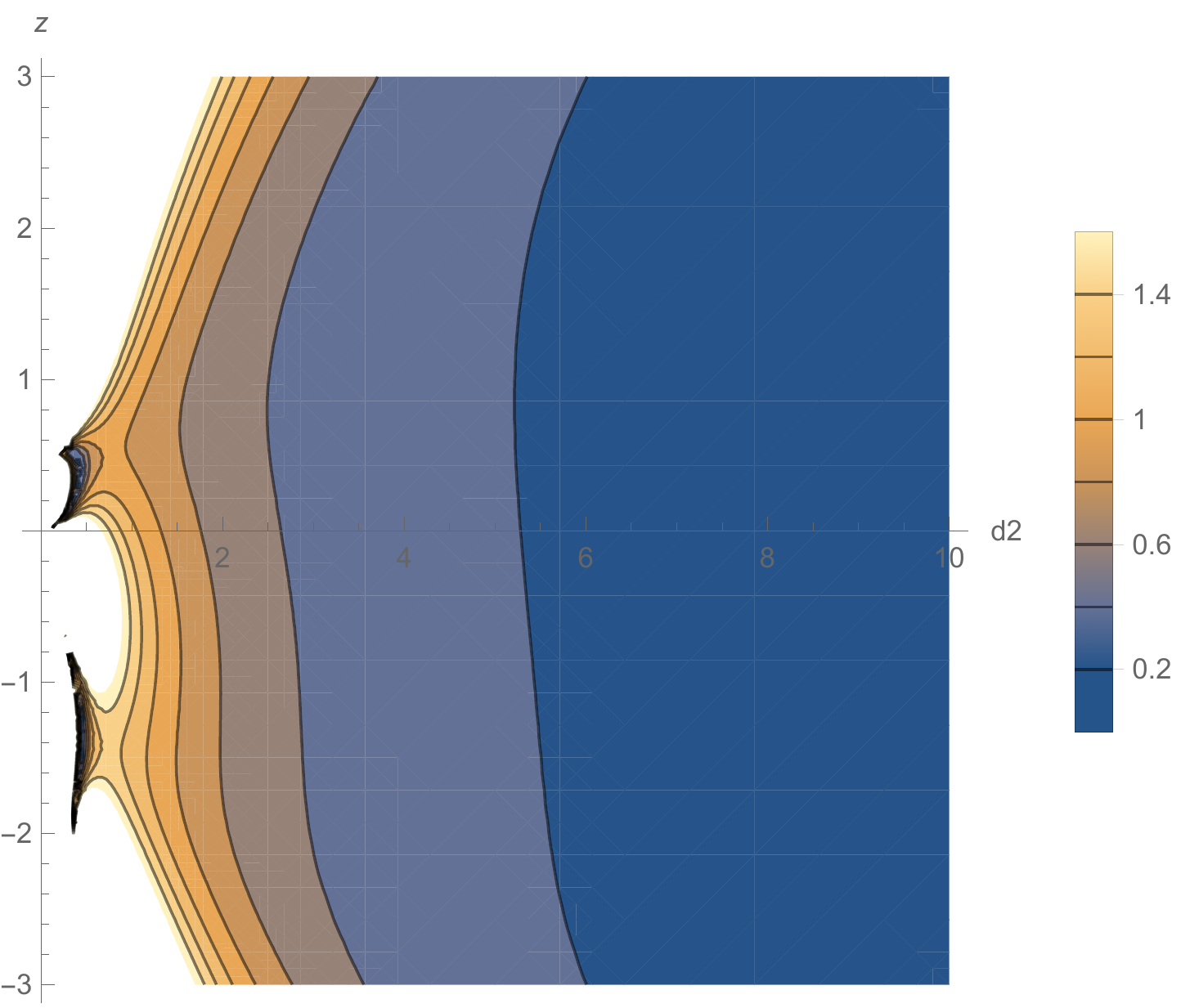}}
      
   \caption{Direct accuracy analysis of the embedded soliton.}\label{Fig5}
\end{figure*}
\begin{figure}[htp]
	\centering
	\includegraphics[width=.5\textwidth]{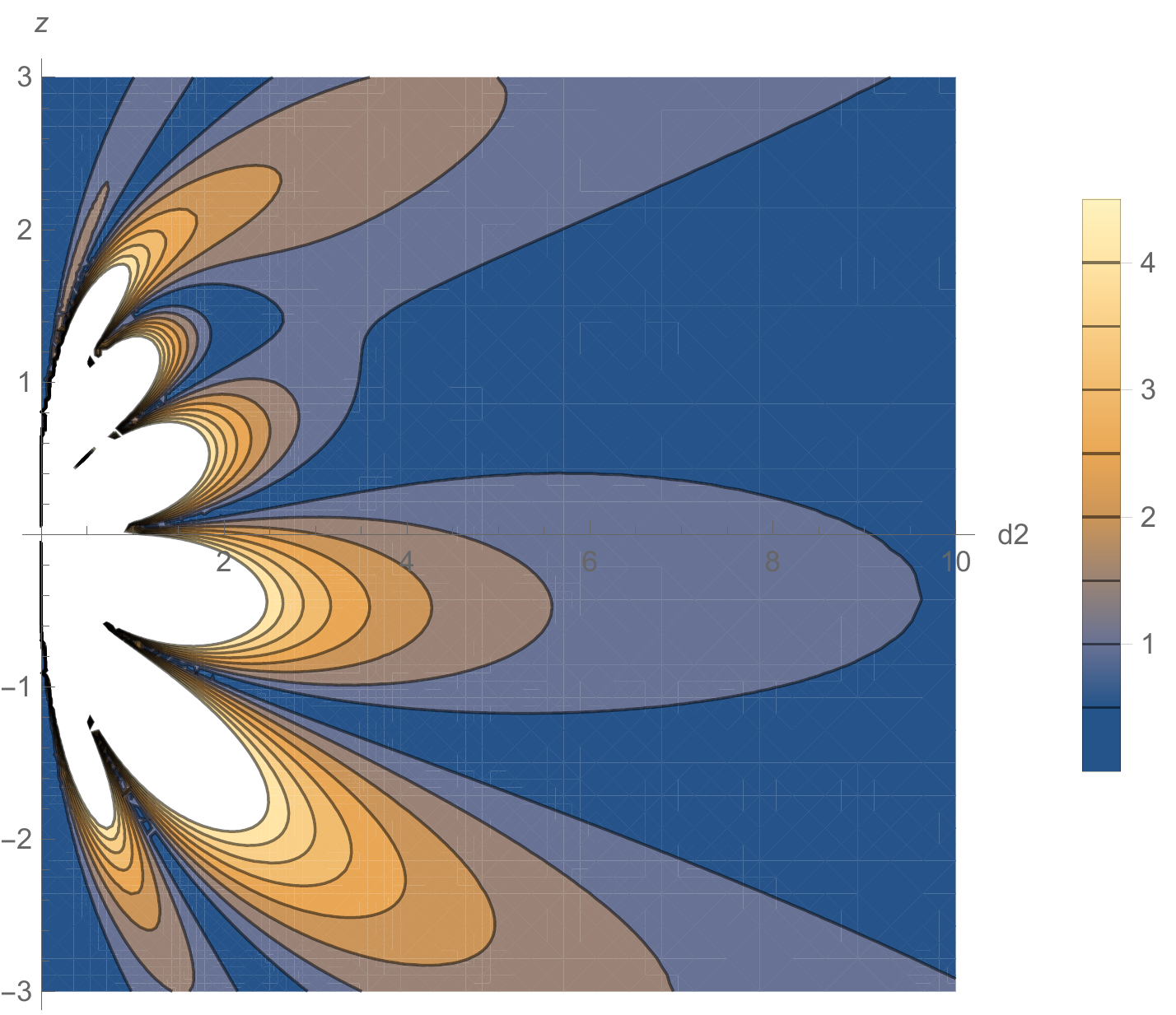}
	\caption{Direct accuracy analysis of the embedded soliton for $\alpha=0$}\label{Fig6}
\end{figure}

In Figures \ref{Fig5} and \ref{Fig6} a direct analysis of the accuracy is carried out for $\mu=1$ and $z$ varying from-3 to 3 while $d_2$ goes from 0 to 10. Figures \ref{Fig5a}, \ref{Fig5b} and \ref{Fig5c} represent the results when $\alpha=1$ and $d_1=-3$, $-1$ and $-1/3$ respectively. Figure \ref{Fig6} presents the results when $\alpha=0$ and $d_1=-1$.

\clearpage


\section{A General Class of Lagrangians and the Associated variational ODEs}
\label{sec:3}

Now, let us generalize the previous ideas to a class of Lagrangians. To do so, consider that the fourth derivative term in a nonlinear ODE to be of the very general form
\begin{equation}
	c_1(u,u')c_2(u'')u^{(4)},
\end{equation}
and we want to match it to the Euler--Lagrange equation
\begin{multline}\label{Eq:EL}
	u^{(4)}L_{u''u''}+\left(u^{(3)}\right)^2L_{u''u''u''}+2u''u^{(3)}L_{u'u''u''}-u''(L_{u'u'}-L_{uu''}) \\
	+\left(u''\right)^2L_{u'u'u''}+2u'u^{(3)}L_{uu''u''}-u'L_{uu'}+2u'u''L_{uu'u''}+\left(u'\right)^2L_{uuu''}+L_{u}=0.
\end{multline}

Thus, at order $O\left(u^{(4)}\right):$ we get
\begin{equation}\label{Eq:Coeff2}
	L_{u''u''}=c_1(u,u')c_2(u''),
\end{equation}
or, on integrating, 
\begin{equation*}
	L_{u''}=c_1(u,u')\int c_2(u'') du''+c_3(u,u').
\end{equation*}
A second integration now yields:
\begin{equation}\label{Eq:Coeff3}
	L=c_1(u,u')c_4(u'')+c_3(u,u')u''+c_5(u,u'),
\end{equation}
where 
\begin{equation}\label{Eq:IntCond}
	c_4(u'')={\displaystyle \iint c_2(u'')du''du'' }  .
\end{equation}

Hence, we must have
\begin{subequations}\label{Eq:Coeff5}
	\begin{align}
		L_{u''u''u''} & = c_1(u,u') \dfrac{dc_{2}(u'')}{du''}, \label{Eq:Coeff5a}\\
		L_{u'u''u''} & = \dfrac{\partial c_{1}(u,u')}{\partial u'}c_2(u''),\label{Eq:Coeff5b}\\
		L_{uu''u''} & = \dfrac{\partial c_{1}(u,u')}{\partial u}c_2(u''), \label{Eq:Coeff5c}\\
		L_{u} & = \dfrac{\partial c_{1}(u,u')}{\partial u}c_4(u'')+\dfrac{\partial c_{3}(u,u')}{\partial u}u''+\dfrac{\partial c_{5}(u,u')}{\partial u}, \label{Eq:Coeff5d}\\
		L_{u'} & =\dfrac{\partial c_{1}(u,u')}{\partial u'}c_4(u'')+\dfrac{\partial c_{3}(u,u')}{\partial u'}u''+\dfrac{\partial c_{5}(u,u')}{\partial u'}, \label{Eq:Coeff5int}\\
		L_{u'u'} & =\dfrac{\partial^2 c_{1}(u,u')}{\partial u'^2}c_4(u'')+\dfrac{\partial^2 c_{3}(u,u')}{\partial u'^2}u''+\dfrac{\partial^2 c_{5}(u,u')}{\partial u'^2}, \label{Eq:Coeff5e}\\
		L_{uu''} & =\dfrac{\partial c_{1}(u,u')}{\partial u}\dfrac{dc_4(u'')}{d u''}+\dfrac{\partial c_{3}(u,u')}{\partial u}, \label{Eq:Coeff5f}\\
		L_{u'u'u''} & =\dfrac{\partial^2 c_{1}(u,u')}{\partial u'^2}\dfrac{dc_4(u'')}{d u''}+\dfrac{\partial^2 c_{3}(u,u')}{\partial u'^2}, \label{Eq:Coeff5g}\\
		L_{uu'} & =\dfrac{\partial^2 c_{1}(u,u')}{\partial uu'}c_4(u'')+\dfrac{\partial^2 c_{3}(u,u')}{\partial uu'}u''+\dfrac{\partial^2 c_{5}(u,u')}{\partial uu'}, \label{Eq:Coeff5h}\\
		L_{uu'u''} & =\dfrac{\partial^2 c_{1}(u,u')}{\partial uu'}\dfrac{dc_4(u'')}{d u''}+\dfrac{\partial^2 c_{3}(u,u')}{\partial uu'}, \label{Eq:Coeff5i}\\
		L_{uuu''} & =\dfrac{\partial^2 c_{1}(u,u')}{\partial u^2}\dfrac{dc_4(u'')}{d u''}+\dfrac{\partial^2 c_{3}(u,u')}{\partial u^2}. \label{Eq:Coeff5j}
	\end{align}
\end{subequations}

Therefore, Equation \eqref{Eq:EL} defines the most general variational fourth order ODE (which may often be the traveling wave equation of a particular nonlinear PDE of interest \cite{Smi09, VPh}) with coefficients given by \eqref{Eq:Coeff2}--\eqref{Eq:Coeff5}, and can now be written, for completely general variable dependence of each function, as
\begin{multline} \label{vareq}
	u^{(4)}c_1c_2 +\left(u^{(3)}\right)^2c_1c_{2u''} +2u''u^{(3)}c_{1u'}c_2 +2u'u^{(3)}c_{1u}c_2 \\ +\left(u''\right)^2\left(c_{1u'u'}c_{4u''}+c_{3u'u'}\right) -u''\left(c_{1u'u'}c_4+c_{3u'u'}u''+c_{5u'u'}-c_{1u}c_{4u''}-c_{3u}\right) \\
	+2u'u''\left( c_{1uu'}c_{4u''}+c_{3uu'} \right) -u'\left( c_{1uu'}c_4+c_{3uu'}u''+c_{5uu'} \right) \\ 
	+\left(u'\right)^2\left( c_{1uu}c_{4u''}+c_{3uu} \right) +c_{1u}c_4+c_{3u}u''+c_{5u}=0 
\end{multline}

In Section V, we shall prove that \eqref{vareq} generalizes the most general currently discussed variational fourth-order variational ODEs, as derived by Fels\cite{Fels}, and which in fact turn out to be special cases of \eqref{vareq}.

\subsection{A Case Worth Mentioning}
\label{subsec:3.1}

Notice that for the very general case 
\begin{equation*}
	c_1(u,u')=\left( \sum_{k=0}^m a_k u^k \right)\left( \sum_{l=0}^n a_l (u')^l \right) ,
\end{equation*}
one may directly obtain the variational ODE \eqref{vareq}. We consider a specific example of this class of Lagrangians in detail in the following section, before returning to the general case and comparing our results to the most general differential geometric criteria obtained earlier for the existence of Lagrangians of fourth-order ODEs.


\section{A member of the family of generalized Lagrangians}
\label{sec:4}

If in the previous section we set 
\begin{align}
	c_1(u,u') & =-d_1+d_2u+d_3u',&  c_2(u'') & =1,\\  c_3(u,u') & =u^2, & c_5(u,u') & =u'^2,
\end{align}
then the Lagrangian reads
\begin{equation}
	L=\left(-d_1+d_2u+d_3u_z\right)u_{zz}^2+u^2 u_{zz}+u_z^2,
\end{equation}
and its Euler--Lagrange equation is
\begin{equation}\label{odeg}
	2u'^2+2(-1+2u)u''+\dfrac{3d_2{u''}^2}{2}+2d_2u'+u'''-d_1u^{(4)}+d_2uu^{(4)}=0.
\end{equation}

In the next subsections we will apply the variational approaches of Section II to obtain solitary wave solutions for the variational ODEs of
this general class of Lagrangians.

\subsection{Regular solitons}
\label{subsec:4.1}
Using the same ansatz \eqref{gtf} and following the process outlined in subsection \ref{subsec:2.1} one gets the action

\begin{equation}
	S=\dfrac{A^2\left(27\sqrt{2}d_1-18\sqrt{2}\rho^2+16\sqrt{3}A\left( -d_2+\rho^2 \right)\right)\sqrt{\pi}}{36\rho^3}.
\end{equation}

\noindent
By solving the system of equations generated by differentiating the action with respect to the parameters one gets the non--trivial solution
\begin{equation}
	A  = \frac{3\sqrt{6}}{64d_2}\tilde{\nu} \qquad \text{and} \qquad \rho^2  = \frac{3}{4}\nu,
\end{equation}
where $\nu=7d_1-2d_2+\sqrt{49d_1^2-36d_1d_2+4d_2^2}$ and $\tilde{\nu}=7d_1+2d_2-\sqrt{49d_1^2-36d_1d_2+4d_2^2}.$

\begin{figure}[ht]
	\centering
	\includegraphics[width=.5\textwidth]{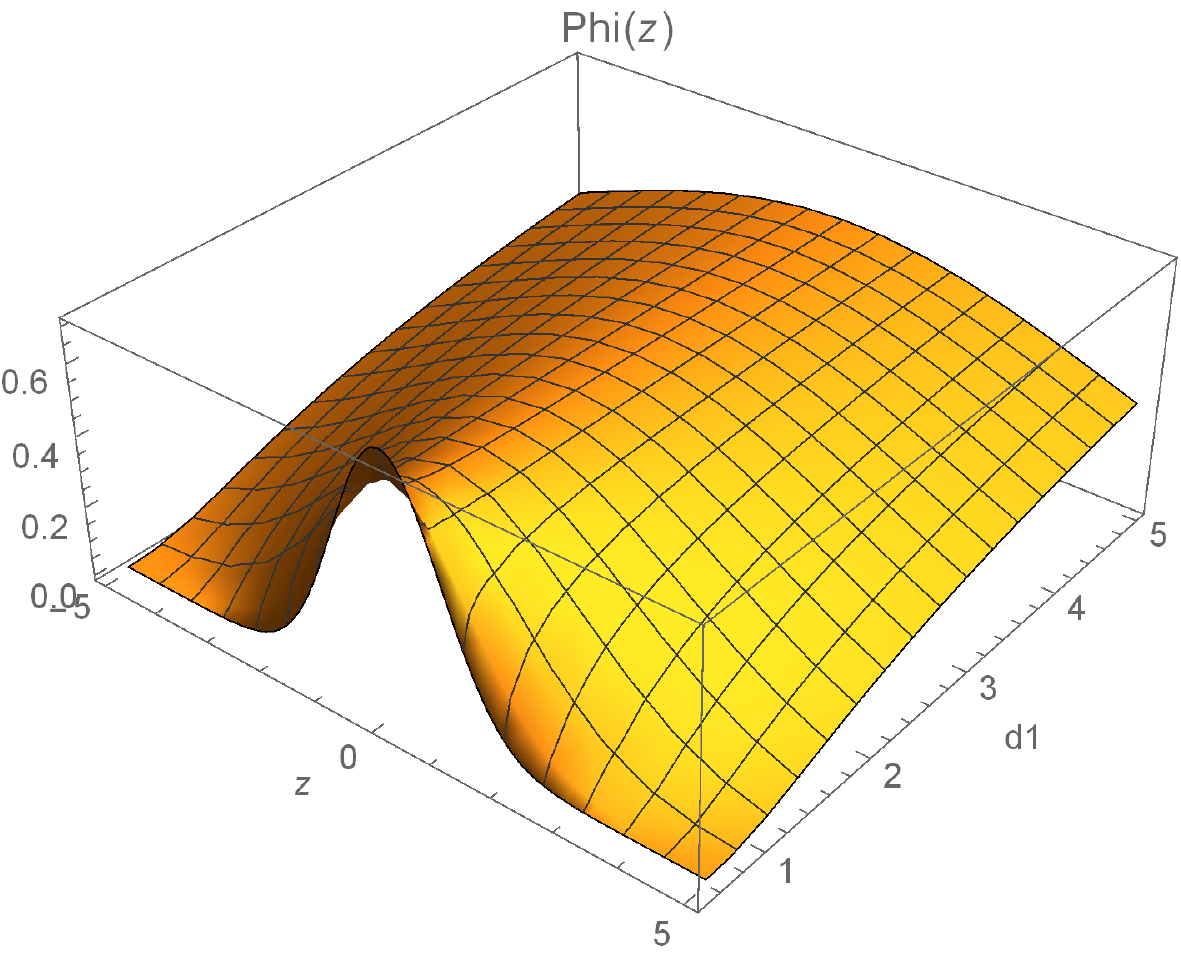}
	\caption{The regular soliton plotted for $d_2=1$.}\label{Fig12}
\end{figure}

In Figure \ref{Fig12} the regular soliton for these values of $\rho$ and $A$ is presented, using $d_2=1$ and varying $z$ from -5 to 5 and $d_1$ from 0.5 to 5.

\subsection{Embedded solitons}
\label{subsec:4.2}
Next, we follow the procedure presented in subsection \ref{subsec:2.2}, we start by considering an ansatz like \eqref{aes}, which yields the averaged action
\begin{multline}
	\frac{\sqrt{\pi } A e^{-\frac{\rho^2}{4 d_1}}}{144 d_1^2 \rho^3}\left(e^{\frac{\rho^2}{4 d_1}}\left( -144\alpha^2(2d_1-d_2)+64\sqrt{3}A^2d_1^2(d_2-\rho^2) \right)\right. \\ \left. -9\sqrt{2}A\left( \alpha \left( -9d_2\rho^4-4d_1^2\left(3d_2-4\rho^2 \right)-12d_1\rho^2 (d_2-2\rho^2) \right) \right.\right. \\
	\left.\left. +4d_1^2 e^{\frac{\rho^2}{4 d_1}}(3d_1-2\rho^2) \right)\right)=S,
\end{multline}

when
\begin{equation}
	\kappa(c)=\frac{\sqrt{2}}{\sqrt{d_1}}.
\end{equation}

\noindent
These yield the non--trivial solutions for $A$ and $\rho^2$:
\begin{equation*}
	A=\frac{\sqrt{\frac{3}{2}}\left(3d_1-\frac{4}{3}\xi\right)}{4\left(d_2-\frac{2}{3}\xi \right)}\quad \text{and} \quad \rho^2=\frac{2}{3}\xi,
\end{equation*}
where $$\xi=\frac{3d_1d_2-4d_1^2-\sqrt{2}\sqrt{d_1^2\left(8d_1^2+24d_1d_2-9d_2^2\right)}}{8d_1-3d_2}.$$

\begin{figure*}[ht]
   \subfloat[The embedded soliton plotted for $\alpha=5$.]{\label{Fig8a}
      \includegraphics[width=.47\textwidth]{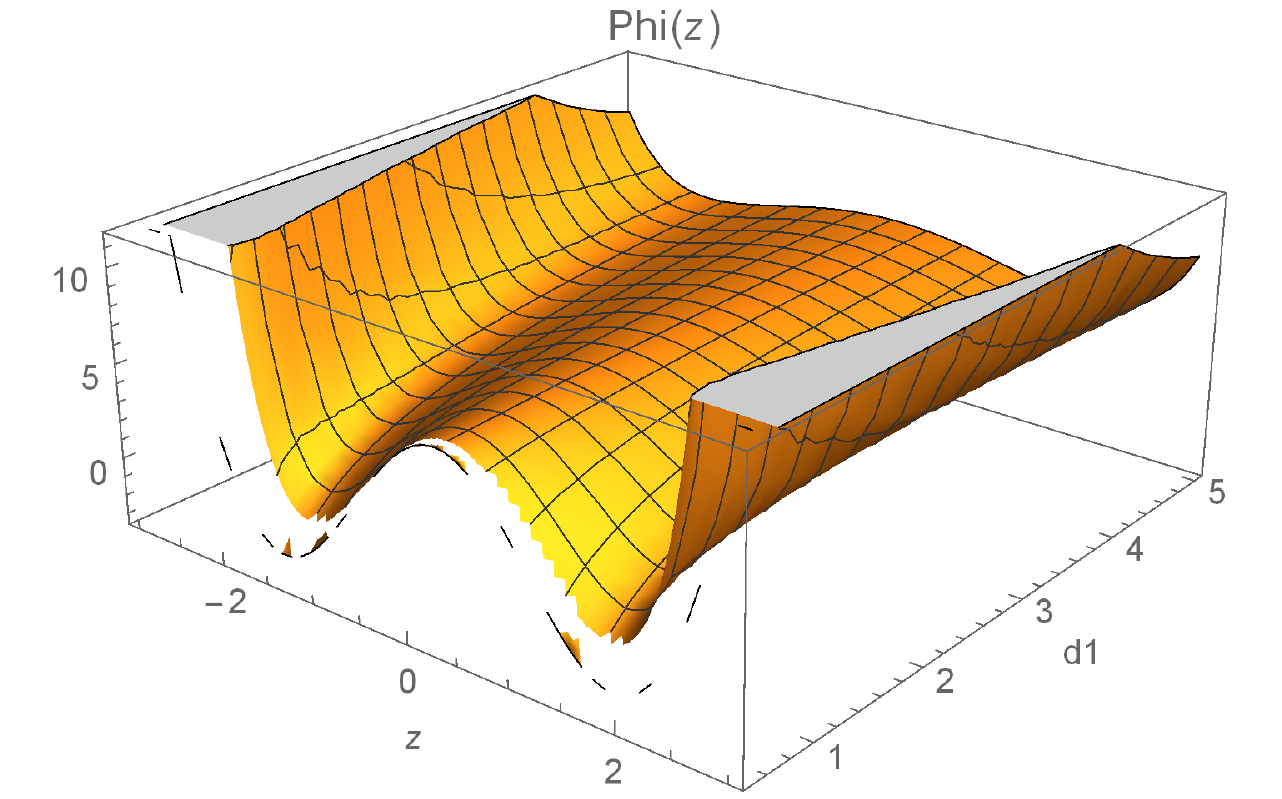}}
~
   \subfloat[The embedded soliton plotted for $d_1=5$.]{\label{Fig8b}
     \includegraphics[width=.47\textwidth]{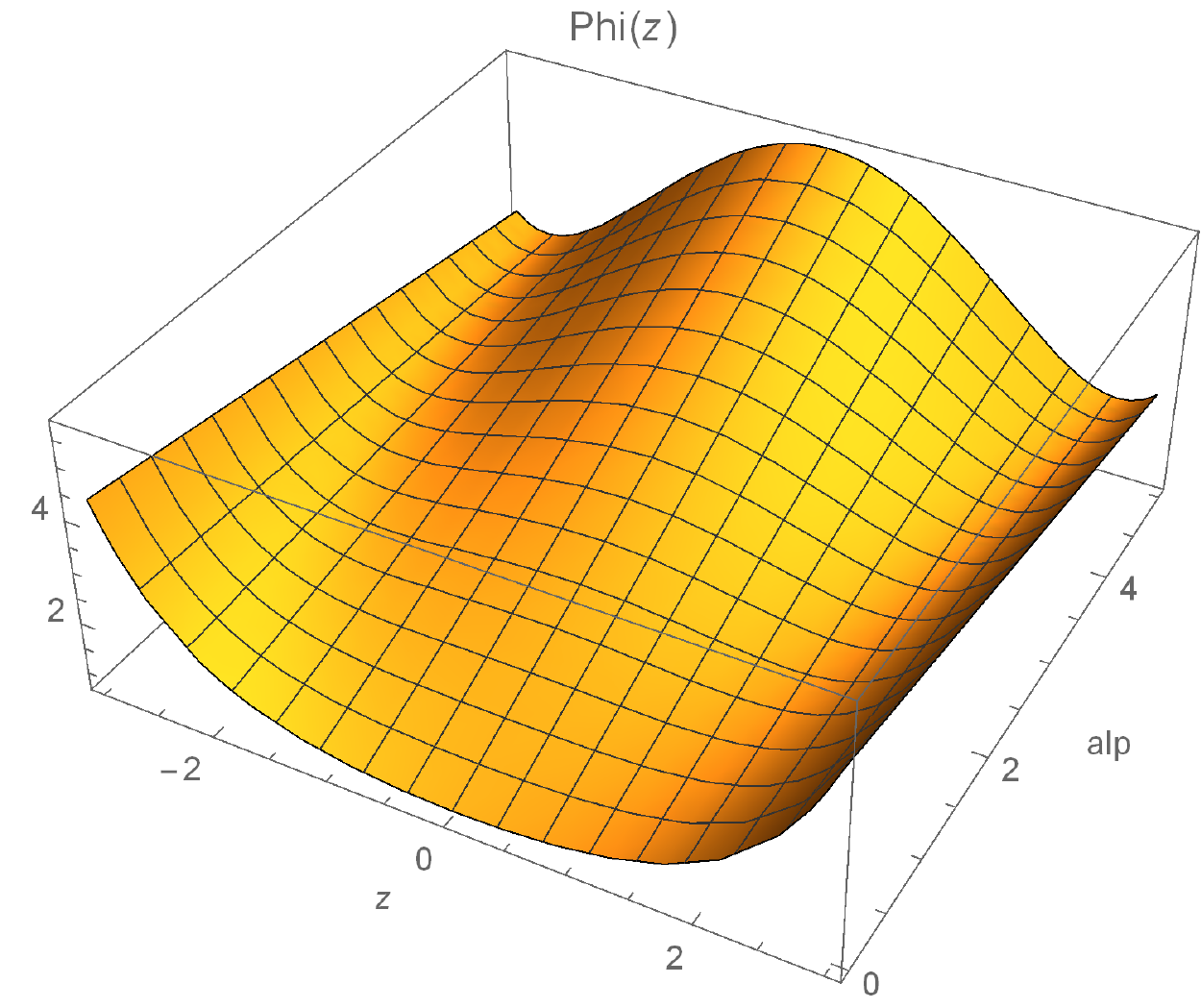}}

   \caption{Embedded soliton for $z$ and different values of the parameters.}\label{Fig8}
\end{figure*}

In Figure \ref{Fig8} the embedded soliton is presented. In it $z$ varies from -3 to 3, with $d_2=2$ and $d_2=5$ respectively. In Figure \ref{Fig8a} $\alpha=5$ while $d_1$ varies between 0.5 to 5. In Figure \ref{Fig8b} $d_1=5$ while $\alpha$ varies from 0 to 5.




\section{Generalizations of previous existence conditions for Lagrangians}
\label{sec:5}

As mentioned in Section 1, the inverse problem for a fourth-order ODE was treated by Fels \cite{Fels}. We now proceed to compare our generalized family of Lagrangians \eqref{Eq:Coeff3}-\eqref{Eq:IntCond} for variational fourth-order ODEs against the criteria obtained by Fels using differential geometric approaches, as well as to various families of variational ODEs based on Fels' criteria \cite{CGK, GC}. 

To this end, we state Fels' principal result here:

\noindent
{\it Theorem: A fourth-order ODE 
	
	\begin{equation}
		u^{(4)} = F(z, u, u', u'')
	\end{equation}
	admits a variational representation with a non-degenerate second order Lagrangian if and only if the following two conditions are satisfied
	
	\begin{equation} \label{Fels1}
		F_{u'''u'''u'''} = 0 
	\end{equation}
	
	and
	
	\begin{equation} \label{Fels2}
		F_{u'} + F_{u'''zz}/2 - F_{u''z} - 3 F_{u'''}F_{u'''z}/ 4 + F_{u''} F_{u'''}/2 + F_{u'''}^3 / 8 = 0
	\end{equation}
	
}

Clearly, our general variational equation \eqref{vareq} corresponds to 

\begin{multline} \label{vareq1}
	F(z, u, u', u") = - ( \left(u^{(3)}\right)^2c_1c_{2u''} +2u''u^{(3)}c_{1u'}c_2 +2u'u^{(3)}c_{1u}c_2 \\ +\left(u''\right)^2\left(c_{1u'u'}c_{4u''}+c_{3u'u'}\right) -u''\left(c_{1u'u'}c_4+c_{3u'u'}u''+c_{5u'u'}-c_{1u}c_{4u''}-c_{3u}\right) \\
	+2u'u''\left( c_{1uu'}c_{4u''}+c_{3uu'} \right) -u'\left( c_{1uu'}c_4+c_{3uu'}u''+c_{5uu'} \right) \\ 
	+\left(u'\right)^2\left( c_{1uu}c_{4u''}+c_{3uu} \right) +c_{1u}c_4+c_{3u}u''+c_{5u}) / (c_1 c_2) 
\end{multline}

It is now straightforward to check that this general $F(z, u, u', u")$ for our generalized family of Lagrangians (which we copy here for ease of understanding of the following discussion) 

\begin{equation}\label{Eq:Coeff3_5}
	L=c_1(u,u')c_4(u'')+c_3(u,u')u''+c_5(u,u'),
\end{equation}
where 
\begin{equation}\label{Eq:IntCond_5}
	c_4(u'')={\displaystyle \iint c_2(u'')du''du'' }  ,
\end{equation}

\noindent
satisfies the first condition \eqref{Fels1} above for arbitrary functions $c_i$ in the Lagrangian \eqref{Eq:Coeff3_5}.

However, the second condition \eqref{Fels2} above is not satisfied for arbitrary functions $c_i$ in the Lagrangian \eqref{Eq:Coeff3_5}. It is however satisfied for the special case
\begin{equation}
	c_1(u,u') = 1, c_2(u'') = 1,
\end{equation}
considered in Fels' \cite{Fels} differential geometric derivation of this criterion, as well as in the treatments of various models based on this condition \eqref{Fels2} \cite{CGK, GC}.

{\it This leads us to the following two conclusions about the two ways that our generalized family of Lagrangians \eqref{Eq:Coeff3_5}-\eqref{Eq:IntCond_5} and the associated variational ODEs \eqref{vareq} are more general than the criteria developed in Fels\cite{Fels}:\\
	
	\noindent
	a. our Lagrangian has four arbitrary or free functions $c_1, c_3, c_4$, and $c_5$, in place of Fels' single function $F$, which, in our case is given by \eqref{vareq1};
	
	\noindent
	and
	
	\noindent
	b. the leading coefficient $c_1(u,u') c_4(u'')$ in our variational ODE \eqref{vareq} may be far more general than for Fels or other models based on his criteria for a variational
	representation.
}


\section{Conclusions}
\label{sec:conclusions}

In conclusion, we have derived a novel, significantly generalized hierarchy of Lagrangians for fourth-order nonlinear ODEs, and employed them to construct families of regular and embedded solitary waves of some member equations.

The method of derivation of the family of Lagrangians is both straightforward and novel. In particular, our family of Lagrangians and the associated variational ODEs are more general than the most general variational formulations derived earlier for fourth-order ODEs \cite{Fels} in two significant ways:

\noindent
a. our Lagrangian has four arbitrary or free functions $c_1, c_3, c_4$, and $c_5$, in place of the  single function $F$ in the earlier general variational criteria 

\noindent
and

\noindent
b. the leading coefficient $c_1(u,u') c_4(u")$ in our variational ODEs \eqref{vareq} may be far more general than in earlier work, or in other models based on the earlier general criteria for a variational representation.

Future work on extending the approach developed here to sixth and eighth-order variational ODE models is ongoing.

\begin{acknowledgements}
The authors would like to thank the SIAM Student Chapter at UCF.
\end{acknowledgements}

%
%


\begin{thebibliography}{}
%
%

\bibitem{PUM} Manheim, P and Davidson, A, Dirac quantization of the Pais-Uhlenbeck fourth order oscillator, Phys Rev A71 (2005) 042110; 

\bibitem{Man} Manheim, P, Solution to the ghost problem in fourth order derivative theories, Foundations Phys 37 (2007) 532.

\bibitem{TCG} Tanriver, U, Roy Choudhury, S and Gambino, G, Lagrangian dynamics and possible isochronous behavior in several classes of non-linear second order oscillators via the use of the Jacobi last multiplier, Intl, J. Nonlin. Mech. 74 (2015) 100.

\bibitem{CGK} Ghosh Choudhury, A, Guha, P and Kudryashov, N, A Lagrangian description of the higher-order Painlev\'e equations, Appl. Math. Comput. 218 (2012) 6612.

\bibitem{VK} Vogel, T and Kaup, DJ, Quantitative measurement of variational approximations, Physics Letters A362 (2007)289.

\bibitem{Helm} Helmholtz, H, Ueber die physikalische bedeutung des prinicips der kleinsten wirkung, J. Reine Angew Math, 100 (1887) 137.

\bibitem{Lop} Lopuszanski, J, The inverse variational problem in classical mechanics (World Scientific, Singapore, 1999).

\bibitem{Dar} Darboux, G, Lecon sur la theorie generale des surfaces (Gauthier-Villars, Paris, 1894).

\bibitem{Fels} Fels, M, The inverse problem of the calculus of variations for scalar fourth-order ordinary differential equations, Trans. Amer. Math Soc, 348 (1996) 5007.

\bibitem{Ili07} Ilison, L, Salupere, A, and Peterson, P, On the propagation of localized perturbations in media with microstructure, Sci. Phys. Math, 56 (2007) 84.

\bibitem{Let09} Leto, J A, and Choudhury, S R, Solitary wave families of a generalized microstructure PDE, Commun. Nonlinear Sci. Numer. Simul, 14 (2009) 1999.

\bibitem{Smi09} Smith, T B, and Choudhury, S R, Regular and embedded solitons in a generalized Pochammer PDE, Commun. Nonlinear Sci. Numer. Simul, 14 (2009) 2637.



\bibitem{VPh} Vogel, T, Soliton Solutions Of Nonlinear Partial Differential Equations Using Variational Approximations And Inverse Scattering Techniques, PhD Thesis, University of Central Florida (2007)

\bibitem{GC} Guha, P and Ghosh Choudhury, A, On Lagrangians and Hamiltonians of some fourth-order nonlinear Kudryashov ODEs, Communications in Nonlinear Science and Numerical Simulation, 16, (2011) 3914.

\end{thebibliography}


\end{document}